\documentclass[11pt,twoside,a4paper]{article}
\usepackage{amsmath,amssymb,amsthm,pifont}

\usepackage{graphicx}
\newtheorem{thm}{Theorem}[section]
\newtheorem{lem}[thm]{Lemma}
\newtheorem{pro}[thm]{Proposition}
\newtheorem{cor}[thm]{Corollary}

\newtheorem{claim}{Claim}

\theoremstyle{definition}

\theoremstyle{remark}
\newtheorem{rem}[thm]{Remark}

\newcommand{\R}{\mathbb{R}}
\newcommand{\Z}{\mathbb{Z}}
\newcommand{\N}{\mathbb{N}}

\newcommand{\cA}{\mathcal{A}}

\newcommand{\cN}{\mathcal{N}}

\newcommand{\cU}{\mathcal{U}}
\newcommand{\cV}{\mathcal{V}}
\newcommand{\cW}{\mathcal{W}}

\newcommand{\al}{\alpha}
\newcommand{\be}{\beta}

\newcommand{\Ga}{\Gamma}
\newcommand{\de}{\delta}
\newcommand{\De}{\Delta}
\newcommand{\ep}{\varepsilon}

\newcommand{\si}{\sigma}

\newcommand{\la}{\lambda}

\renewcommand{\phi}{\varphi}
\newcommand{\lo}[1]{\ensuremath{\overline{#1}}}

\newcommand{\dist}{\operatorname{dist}}
\newcommand{\diam}{\operatorname{diam}}

\newcommand{\hyp}{\operatorname{H}}
\newcommand{\id}{\operatorname{id}}
\newcommand{\Lip}{\operatorname{Lip}}

\newcommand{\lasdim}{\operatorname{\ell-asdim}}
\newcommand{\asdim}{\operatorname{asdim}}
\newcommand{\hypdim}{\operatorname{hypdim}}
\newcommand{\ba}{\operatorname{ba}}
\newcommand{\st}{\operatorname{st}}
\newcommand{\cone}{\operatorname{Co}}
\newcommand{\pt}{\operatorname{pt}}
\newcommand{\mesh}{\operatorname{mesh}}

\newcommand{\an}{\operatorname{An}}
\newcommand{\ldim}{\operatorname{\ell-dim}}

\newcommand{\es}{\emptyset}
\renewcommand{\d}{\partial}
\newcommand{\di}{\d_{\infty}}
\newcommand{\set}[2]{\{#1:\,\text{#2}\}}
\newcommand{\sm}{\setminus}
\newcommand{\sub}{\subset}
\newcommand{\sups}{\supset}

\newcommand{\ov}{\overline}
\newcommand{\wt}{\widetilde}
\newcommand{\wh}{\widehat}
\newcommand{\md}{\!\mod}

\begin{document}

\title{Dimensions of products of hyperbolic spaces}
\author{by Nina Lebedeva {\footnote{supported by RFFI Grant
05-01-00939}}}

\date{}
\maketitle

\begin{abstract}
We give estimates on asymptotic dimensions of
 products of general hyperbolic spaces with following
applications to the hyperbolic groups.
We give examples of strict inequality in the product theorem
for the asymptotic dimension in the class of the hyperbolic groups;
and examples of strict inequality in the product theorem
for the hyperbolic dimension.
We prove that $\R$ is dimensionally full for the
asymptotic dimension in the class of the hyperbolic groups.

\end{abstract}

\section{Introduction}
The purpose of this article is to answer
some questions in asymptotic dimension theory (in
particular posed in \cite{Dr3}).
For this we give different estimates on asymptotic dimensions of
general
hyperbolic spaces and products of hyperbolic spaces and apply them
to hyperbolic groups.

To formulate our results we need some definitions.
We say that metric spaces
$\set{X_\al}{$\al\in\cA$}$
have the linearly controlled dimension
$\ldim\le k$
uniformly if there is a constant
$\de\in(0,1)$
such that for every sufficiently small
$\tau>0$
and for every
$\al\in\cA$
there exists a
$(k+1)$-colored
open covering
$\cU_\al$
of
$X_\al$
with
$\mesh(\cU_\al)\le\tau$
and with Lebesgue number
$L(\cU_\al)\ge\de\tau$.

For a metric space
$X$
and
$a>0$,
we denote by
$aX$
the metric space obtained from
$X$
by multiplying all distances by
$a$.

In sect.~\ref{subsect:defcapdim}, we recall the definitions
of the asymptotic dimension
$\asdim$
and the linearly controlled asymptotic dimension
$\ell$-$\asdim$
used in the following theorem, these are some quasi-isometry
invariants of metric spaces.

\begin{thm}\label{thm:asprod<} Let
$X_1$, $X_2$
be visual Gromov hyperbolic spaces and assume that
the metric spaces of the family
$\set{a\di X_1\times b\di X_2}{$a\ge 1, b\ge 1$}$
have
$\ldim\le k$
uniformly, where the boundaries at infinity
$\di X_1$, $\di X_2$
are taken with some visual metrics. Then
$$\asdim(X_1\times X_2)\le k+2.$$
\end{thm}

The property of a Gromov hyperbolic space
$X$
to be visual is a rough version of the property that every
point
$x\in X$
lies on a geodesic ray emanating from a fixed point
$x_0\in X$;
for the precise definition see sect.~\ref{sect:pre}.

\begin{rem} A similar result holds true in the case of
an arbitrary finite number
$\ge 2$
of factors, and furthermore the estimate remains true if we replace
$\asdim$
by
$\lasdim$,
that is:

Let
$X_1,\ldots,X_n$
be visual Gromov hyperbolic spaces and assume that
the metric spaces of the family
$\set{a_1\di X_1\times\cdots\times a_n\di X_n}{$a_i\ge 1, i\in\{1,\ldots,n\}$}$
have
$\ldim\le k$
uniformly. Then
$$\asdim(X_1\times\cdots\times X_n)\le k+n$$
and
$$\lasdim(X_1\times\cdots\times X_n)\le k+n.$$

For simplicity of notations, we shall consider the case of
two spaces and the asymptotic dimension
$\asdim$.
\end{rem}

\begin{thm}\label{thm:asprod>} Let
$X_1,\ldots,X_n$
be geodesic Gromov hyperbolic spaces. Then
$$\asdim(X_1\times\cdots\times X_n)
  \ge\ldim(\di X_1\times\cdots\times\di X_n\times[0,1]^n).$$
\end{thm}

\begin{cor}\label{cor:G} Let
$\Ga_1$, $\Ga_2$
be  Gromov hyperbolic groups. Then
$$\asdim(\Ga_1\times\Ga_2 )=\dim(\di\Ga_1\times\di\Ga_2)+2.$$
\end{cor}

As an application of  this corollary, we

1) give examples of strict inequality in the product theorem
for the asymptotic dimension in the class of hyperbolic groups
(cf. \cite[Problem 23]{Dr3}; a corresponding example
in the class of general metric spaces is given in
\cite{BL});

2) give examples of strict inequality in the product theorem
for the hyperbolic dimension;

3) prove the equality
$\asdim(\Ga\times\R)=\asdim(\Ga)+1$
for any hyperbolic group
$\Ga$.

The last equality is not true for general metric spaces.
An example of a metric space
$X$
of bounded geometry with finite asymptotic dimension
for which
$\asdim(X\times\R)=\asdim X$
is constructed in \cite{Dr4}, and the question is,
if the equality above is true for groups
(\cite[Problem 21]{Dr3}).

\begin{thm}\label{thm:lasdim} Let
$X$
be a visual Gromov hyperbolic space. Then
$$\lasdim X\le\ldim(\di X\times[0,1]).$$
\end{thm}

This estimate strengthens the estimate
$\asdim X\le\ldim\di X+1$
proved in \cite{Bu1, Bu2}. Unfortunately, we do not
know whether the space
$[0,1]$
is dimensionally full for the linearly controlled dimension,
i.e. whether
$\ldim(Y\times[0,1])=\ldim Y+1$
for any (compact or proper) metric space
$Y$.

\begin{cor}\label{cor:lasdim} Let
$X$
be a visual, geodesic Gromov hyperbolic space. Then
$$\lasdim X=\asdim X=\ldim(\di X\times[0,1]).$$
\end{cor}

The first equality does not hold in general and the question is if it
holds for all finitely presented groups (\cite[Problem 41]{Dr3}).

\noindent {\bf  Aknowledgement}
I would like   to thank  Sergei Buyalo for 
many useful remarks and attention to this paper.
I'm also
pleased to acknowledge Max-Plank Mathematical Institute 
and University of Muenster for support and excellent
working conditions while writing   the paper.

\section{Preliminaries}\label{sect:pre}

Here, we recall notions and facts necessary for the paper.

\subsection{Coverings}\label{subsect:coverings}

Let
$Z$
be a metric space. For
$U$, $U'\sub Z$
we denote by
$\dist(U,U')$
the distance between
$U$
and
$U'$,
$\dist(U,U')=\inf\set{|uu'|}{$u\in U,\ u'\in U'$}$
where
$|uu'|$
is the distance between
$u$, $u'$.
For
$r>0$
we denote by
$B_r(U)$
the open
$r$-neighborhood
of
$U$, $B_r(U)=\set{z\in Z}{$\dist(z,U)<r$}$,
and by
$\ov B_r(U)$
the closed
$r$-neighborhood
of
$U$, $\ov B_r(U)=\set{z\in Z}{$\dist(z,U)\le r$}$.
We extend these notations over all real
$r$
putting
$B_r(U)=U$
for
$r=0$,
and defining
$B_r(U)$
for
$r<0$
as the complement of the closed
$|r|$-neighborhood
of
$Z\sm U$,
$B_r(U)=Z\sm\ov B_{|r|}(Z\sm U)$.

Given a family
$\cU$
of subsets in a metric space
$Z$
we define
$\mesh(\cU)=\sup\set{\diam U}{$U\in\cU$}$.
The {\em multiplicity} of
$\cU$, $m(\cU)$,
is the maximal number of members of
$\cU$
with nonempty intersection. We say that a family
$\cU$
is {\em disjoint} if
$m(\cU)=1$.

A family
$\cU$
is called a {\em covering} of
$Z$
if
$\cup\set{U}{$U\in\cU$}=Z$.
A covering
$\cU$
is said to be {\em colored} if it is the union
of
$m\ge 1$
disjoint families,
$\cU=\cup_{a\in A}\cU^a$, $|A|=m$.
In this case we also say that
$\cU$
is
$m$-colored.
Clearly, the multiplicity of a
$m$-colored
covering is at most
$m$.

Let
$\cU$
be a family of open subsets in a metric space
$Z$
which cover
$A\sub Z$.
Given
$z\in A$,
we let
$$L(\cU,z)=\sup\set{\dist(z,Z\sm U)}{$U\in\cU$}$$
be the Lebesgue number of
$\cU$
at
$z$, $L(\cU)=\inf_{z\in A}L(\cU,z)$
the Lebesgue number of the covering
$\cU$
of
$A$.
For every
$z\in A$,
the open ball
$B_r(z)\sub Z$
of radius
$r=L(\cU)$
centered at
$z$
is contained in some member of the covering
$\cU$.

We shall use the following obvious fact (see e.g. \cite{Bu1}).

\begin{lem}\label{lem:insidecov} Let
$\cU$
be an open covering of
$A\sub Z$
with
$L(\cU)>0$.
Then for every
$s\in(0,L(\cU))$
the family
$\cU_{-s}=B_{-s}(\cU)$
is still an open covering of
$A$.
\qed
\end{lem}

\subsection{Hyperbolic spaces}\label{subsect:hypspace}

Here we recall necessary facts from the hyperbolic spaces theory.
For more details one can see e.g. \cite{BoS}.
We also assume that the reader is familiar with notions
like of a geodesic metric space, a geodesic ray etc.

Let
$X$
be a metric space. We use notation
$|xx'|$
for the distance between
$x$, $x'\in X$. For
$o\in X$
and for
$x$, $x'\in X$,
put
$(x|x')_o=\frac{1}{2}(|xo|+|x'o|-|xx'|)$.
The number
$(x|x')_o$
is nonnegative by the triangle inequality, and it is
called the Gromov product of
$x$, $x'$
w.r.t.
$o$.

A  metric space
$X$
is called
$\de$-{\em hyperbolic},
$\de\ge 0$,
if the following inequality
$$(x|y)_o\ge\min\{(x|z)_o,(z|y)_o\}-\de$$
holds for every base point
$o\in X$
and all
$x$, $y$, $z\in X$.

Let
$X$
be a hyperbolic space, i.e.
$X$
is
$\de$-hyperbolic for some
$\de\ge 0$.
We denote by
$\di X$
the (Gromov) boundary of
$X$
at infinity. The Gromov product based at
$o\in X$
naturally extends to
$\di X$.

If a hyperbolic space
$X$
is  geodesic then every geodesic ray in
$X$
represents a point at infinity. Conversely, if
a geodesic hyperbolic space
$X$
is proper (i.e. closed balls in
$X$
are compact), then
every point at infinity is represented by a
geodesic ray.

A metric
$d$
on the boundary at infinity
$\di X$
of
$X$
is said to be {\em visual}, if there are
$o\in X$, $a>1$
and positive constants
$c_1$, $c_2$,
such that
$$c_1a^{-(\xi|\xi')_o}\le d(\xi,\xi')\le c_2a^{-(\xi|\xi')_o}$$
for all
$\xi$, $\xi'\in\di X$.
In this case, we say that
$d$
is visual with respect to the base point
$o$
and the parameter
$a$.
The boundary at infinity is bounded and complete w.r.t.
any visual metric, moreover, if
$X$
is proper then
$\di X$
is compact. If
$a>1$
is sufficiently close to 1, then a visual metric with
respect to
$a$
does exist.

A hyperbolic space
$X$
is called {\em visual}, if for some base point
$x_0\in X$
there is a positive constant
$D$
such that for every
$x\in X$
there is
$\xi\in\di X$
with
$|xx_0|\le(x|\xi)_{x_0}+D$
(one easily sees that this property is independent of
the choice of
$x_0$).

\subsection{Definitions of linearly controlled and
asymptotic dimensions}\label{subsect:defcapdim}

We shall say that a covering
$\cU$
is an
$(L,M,n)$-covering if it is
$n$-colored
with
$\mesh(\cU)\le M$
and
$L(\cU)\ge L$.

The notion of the linearly controlled metric dimension
of a metric space is introduced in \cite{Bu1} where
it is called the `capacity dimension'. One of a number
equivalent definitions of this notion is the following.

The linearly controlled dimension or
$\ell$-dimension
of a metric space
$Z$, $\ldim(Z)$,
is the minimal integer
$m\ge 0$
with the following property: There is a constant
$\de\in(0,1)$
such that for every sufficiently small
$\tau>0$
there exists an open
$(\de\tau, \tau, m+1)$-covering of $Z$.
Here (and in all similar definitions of a dimension below)
we assume that
$\ldim Z=\infty$,
if the number
$m$
with the indicated property does not exist.

The asymptotic dimension is a quasi-isometry invariant of
a metric space introduced in \cite{Gr}.
There are also several equivalent definitions,
see \cite{Gr}, \cite{BD}, and we use the following one.
The asymptotic dimension of a metric space
$X$
is the minimal integer
$\asdim X=m$
such that for every positive
$d$
there is an open covering
$\cU$
of
$X$
with
$m(\cU)\le m+1$, $\mesh(\cU)<\infty$
and
$L(\cU)\ge d$.

The linearly controlled asymptotic dimension is
also a quasi-isometry invariant of a metric space.
We use the following definition.
The linearly controlled asymptotic dimension of a metric space
$X$
is the minimal integer
$\lasdim X=n$
with the following property: There is a constant
$C>1$
such that for every sufficiently large
$d$
there is an open
$(d,Cd,n+1)$-covering $\cU$
of
$X$.

\subsection{Coverings of the product}
Throughout the paper, we assume that the product of metric spaces
$X\times Y$
is endowed with
$l_\infty$-metric,
that is
$|(x,y)(x',y')|=\max\{|xx'|,|yy'|\}$
for each
$(x,y)$, $(x',y')\in X\times Y$.

Let
$\cU$
be an open covering of the product
$X\times Y$
of metric spaces. We say that a pair
$(p,q)$, $p$, $q\ge 0$,
dominates
$\mesh^\times(\cU)$
and write
$\mesh^\times(\cU)\le(p,q)$
or
$(p,q)=\mesh^\boxtimes(\cU)$,
if for each member of
$\cU$
there exists a point
$(x,y)\in X\times Y$
so that this member is contained in the product
$\ov B_p(x)\times\ov B_q(y)$.

We also say that the function
$f(z)\in[0,\infty)\times[0,\infty)$, $z\in X\times Y$,
dominates the pointed
$\mesh^\times(\cU,z)$
and write
$\mesh^\times(\cU,z)\le f(z)$
or
$f(z)=\mesh^\boxtimes(\cU,z)$,
if for each member of
$\cU$
that contains
$z$
there exists a point
$(x,y)\in X\times Y$
so that this member is contained in
$\ov B_p(x)\times\ov B_q(y)$,
where
$(p,q)=f(z)$.
In this sense, the
$\mesh^\times$
(the pointed
$\mesh^\times$)
is determined basically by the meshes
(the meshes at the projected points) of the coverings
projected to the factors.

Similarly, we say that the {\em Lebesgue number}
$L^\times$
at a point
$z=(x,y)\in X\times Y$
of the covering
$\cU$
{\em dominates}
$f(z)=(p,q)\in[0,\infty)\times[0,\infty)$,
$L^\times(\cU,z)\ge f(z)$,
if the product of the balls
$\ov B_p(x)\times\ov B_q(y)$
is contained in some member of the covering.
In this case, we also write
$f(z)=L^\boxtimes(\cU,z)$.

We write
$L^\times(\cU)\ge(p,q)$
or
$(p,q)=L^\boxtimes(\cU)$,
if
$L^\times(\cU,z)\ge(p,q)$
for every
$z\in X\times Y$.
There is no way in general to recover
$L^\times$
from the Lebesgue numbers of the projected coverings.

For the product
$X\times Y$,
it is straightforward to check that the following
two properties are equivalent.

(1) The spaces of the family
$\set{aX\times bY}{$a\ge1, b\ge1$}$
have
$\ldim\le k$
uniformly.

(2) There is a constant
$\de\in(0,1)$
such that for every sufficiently small
$\tau>0$
and for every
$\al$, $\be\in(0,1]$
there exists a
$(k+1)$-colored
open covering
$\cU$
of
$ X\times  Y$
with
$\mesh^\times(\cU)\le(\al\tau,\be\tau)$
and
$L^\times(\cU)\ge(\al\de\tau,\be\de\tau)$.

Let
$\cU$, $\cV$
be open coverings of
$X\times Y$.
We write
$\mesh^\times(\cU)\le(p,q)L^\times(\cV)$
for some
$p$, $q\ge 0$,
if
$\mesh^\times(\cU)\le(pr_1,qr_2)$
for some
$(r_1,r_2)=L^\boxtimes(\cV)$.
In the case
$p=q$,
we write
$\mesh^\times(\cU)\le pL^\times(\cV)$.
Note that if
$\mesh^\times(\cU)\le L^\times(\cV)$
then the covering
$\cU$
is {\em inscribed} in the covering
$\cV$,
that is, every member of
$\cU$
is contained in some member of
$\cV$.

We write
$\mesh^\times(\cU,z)\le\mesh^\times(\cV,z)$,
if every pair
$\mesh^\boxtimes(\cV,z)$
is also a pair
$\mesh^\boxtimes(\cU,z)$.
We write
$L^\times(\cU,z)\ge s L^\times(\cV,z)$
for
$s\ge 0$,
if for any
$p$, $q\ge 0$
with
$L^\times(\cV,z)\ge (p,q)$
we have
$L^\times(\cU,z)\ge (sp,sq)$.

For
$A\sub Z=X\times Y$
and
$r_1$, $r_2>0$,
we denote by
$B_{(r_1,r_2)}(A)$
the union of all products
$B_{r_1}(x)\times B_{r_2}(y)$
with
$(x,y)\in A$,
$B_{(r_1,r_2)}(A)=\cup\set{B_{r_1}(x)\times B_{r_2}(y)}
   {$(x,y)\in A$}$,
and by
$B_{(-r_1,-r_2)}(A)=Z\sm \overline{B}_{(r_1, r_2)}(Z\sm A).$

We also write
$(p,q)\ge (p',q')$
for reals
$p$, $q$, $p'$, $q'$,
if and only if
$p\ge p'$
and
$q\ge q'$.

\section{Auxiliary facts}

\subsection{Separate and qualified families of coverings}

In this section we give a  generalization of
the following lemma, see e.g. \cite{BL}.
\begin{lem}\label{lem:union} Suppose, that
$Z$
is a metric space and
$A$, $B\sub Z$.
Let
$\cU$
be an open covering of
$A$,
$\cV$
an open covering of
$B$
both with multiplicity at most
$m$.
If
$\mesh(\cV)\le L(\cU)/2$
then there exists an open covering
$\cW$
of
$A\cup B$
with multiplicity at most
$m$
and
$\mesh(\cW)\le\max\{\mesh(\cV),\mesh(\cU)\}$,
$L(\cW)\ge\min\{L(\cU)/2,L(\cV)\}$.
\qed
\end{lem}

Let
$\cU_1$, $\cU_2$
be families of subsets in a metric space
$Z$.
We say that
$\cU_1$
and
$\cU_2$
are separate, if every member of
$\cU_1$
intersects no member of
$\cU_2$.

Let
$Z=X\times Y$
be the product of metric spaces and
$Z_{\al}\sub Z$, $\al\in\cA$.
Let
$\cU_{\al}$
be a covering of
$Z_{\al}$.
Let
$S=\{0,\dots,N\}$
be the set of `scales' and
$i:\cA\to S$
a {\em scale function}.

We say that the family
$\cU_{\al}$, $\al\in\cA$
is {\em separate with the scale function}
$i$
if for every
$U\in\cU_{\al}$, $V\in\cU_{\al'}$
with
$\al\neq\al'$
and
$i(\al)=i(\al')$,
we have
$U\cap V=\es$.

We say that the family
$\cU_{\al}$
is {\em qualified by the scale function}
$i$
if the condition
$U\cap V\neq\es$
for some
$U\in\cU_{\al}$, $V\in\cU_{\al'}$
with
$i(\al)<i(\al')$
implies
$L^\boxtimes(\cU_{\al'})\ge 4\mesh^\boxtimes(\cU_{\al})$
for appropriate dominating numbers fixed
for each covering
$\cU_\al$.

\begin{pro}\label{prop:glue} Suppose that
$Z=X\times Y$
is the product of metric spaces and
$Z_{\al}\sub Z$, $\al\in\cA$.
Let
$S=\{0,\dots,N\}$
and let
$i:\cA\to S$
be a scale function. Let
$\cU_{\al}$
be an open
$m$-colored
covering of
$Z_{\al}$
for every
$\al\in\cA$
so that the family of the coverings
$\cU_{\al}$, $\al\in\cA$
is separate with and qualified by the scale function
$i$.
Then there exists an open
$m$-colored
covering
$\cW$
of
$\cup_{\al}Z_{\al}$
with
$\mesh^\times(\cW,z)\le\mesh^\times(\cup_\al \cU_{\al}, z)$,
$L^\times(\cW, z)\ge
L^\times(\cup_\al \cU_{\al}, z)/2$
for every
$z\in\cup_{\alpha}Z_{\al}$.
\end{pro}

\begin{proof} Let
$A$
be the set of colors,
$|A|=m$,
that we may think to be common for all coverings
$\cU_{\al}$,
so that
$\cU_{\al}$
is the union of disjoint families,
$\cU_{\al}=\bigcup_{a\in A}\cU_{\al}^a$.
For
$s\in S$,
we denote by $\cA_s:=\set{\al\in\cA}{$i(\al)=s$}$.

For every
$\al\in\cA$, $a\in A$,
we consider the family
$\cV_\al^a=B_{-L^\boxtimes(\cU_\al)/2}(\cU_\al^a)$,
where the pair
$L^\boxtimes(\cU_\al)$
is dominated by
$L^\times(\cU_\al)$
as in the definition of coverings qualified
by a scale function.

Assume that
$U\in\cU_\al$, $V=B_{-L^\boxtimes(\cU_\al)/2}(V')$
for
$V'\in \cU_{\al'}^a$.
We have

\begin{itemize}
\item[($\ast$)] if
$i(\al)<i(\al')$
and
$U\cap V\neq\es$,
then
$U\cup V\subset V'$.
\end{itemize}
This easily follows from the fact that
$L^\boxtimes(\cU_{\al'})\ge 4\mesh^\boxtimes(\cU_{\al})$
in this case. In particular,
\begin{itemize}
\item[($\ast\ast$)] for every
$U\in\cU_{\al}$,
every
$\al'\in\cA$
with
$i(\al)<i(\al')$
and every
$a\in A$
there  exists at most one
$V\in\cV_{\alpha'}^a$
with
$U\cap V\neq\es$,
\end{itemize}
because the family
$\cU_{\al'}^a$
is disjoint.

We construct an $m$-colored covering $\cW$ by induction over $j\in S$,
independently for each color  $a\in A$.

Fix a color
$a\in A$
and set
$\cW_0^a:=\cup_{\al\in\cA_0}U_\al^a$.
Then the family
$\cW_0^a$
is disjoint and moreover, it possesses properties (2), (3)
below. Assume that for some
$j\in S$
we have already constructed a family
$\cW_j^a$
with the following properties:

\begin{itemize}
\item[(1)] $\cW_j^a$ is disjoint;
\item[(2)] for every
$U\in \cV_\al^a$
with
$i(\al)\le j$
there exists
$W\in\cW^a_j$
so that
$U\subset W$;
\item[(3)] for every
$W\in\cW^a_j$
there  exists
$U\in\cU_\al^a$
with
$i(\al)\le j$
so that
$W\sub U$.
\end{itemize}
In view of ($\ast\ast$), property~(3) implies that for every
$W\in\cW^a_j$
there  exists at most one member
$V\in\bigcup_{\cA_{j+1}}\cV^a_\al$
so that
$V\cap W\neq\es$.
We set
$I(W)=V$
if such
$V$
exists and
$I(W)=\es$
otherwise; so we have the function
$I:\cW^a_j\to\bigcup_{\cA_{j+1}}\cV^a_\al\cup \{\es\}$.

We put
$$\cW^a_{j+1}=\set{W\in\cW^a_j}{$I(W)=\es$}\cup
 \set{(V\cup I^{-1}(V)}{$V\in\bigcup_{\cA_{j+1}}\cV^a_\al$}.$$
The family
$\bigcup_{\cA_{j+1}}\cV^a_\al$
is disjoint. Thus in view of
$(\ast)$
and (1), the family
$\cW^a_{j+1}$
is also disjoint. Furthermore,
$\cW^a_{j+1}$
possesses property (2) (with
$j$
replaced by
$j+1$)
by construction. Property (3) also holds for
$\cW^a_{j+1}$
by the inductive assumption and in view of
($\ast$).

Then proceeding by induction we obtain the family
$\cW^a_{N}$
for each
$a\in A$.
We set
$\cW^a=\cW^a_{N}$
and
$\cW=\cup_a\cW^a$.
Property (2) for
$\cW^a_{N}$
implies that
$\cW$
is a covering of
$\cup_{\al}Z_{\al}$.
By property~(1) for
$\cW^a_{N}$,
$\cW$
is
$m$-colored.
Property~(3) implies
$\mesh^\times(\cW,z)\le\mesh^\times(\cup_\al \cU_{\al},z)$
and property~(2) implies
$L^\times(\cW, z)\ge L^\times(\cup_\al \cU_{\al}, z)/2$
for every
$z\in\cup_{\alpha}Z_{\al}$.
That is,
$\cW$
possesses all required properties.
\end{proof}

\subsection{Barycentric triangulation of products}
\label{subsect:barycentrian}

Recall some standard constructions related to simplicial
polyhedra. Given an index set
$J$,
we let
$R^J$
be the Euclidean space of functions
$J\to\R$
with finite support. For
$x$, $x'\in\R^J$,
the distance
$|xx'|$
is well defined by

$$|xx'|^2 =\sum_{j\in J} (x_j - x_j')^2.$$

Let
$\De^J\sub R^J$
be the standard simplex, i.e.,
$x\in\De^J$
iff
$x_j\ge 0$
for all
$j\in J$
and
$\sum_{\in J} x_j=1$.
The metric of
$R^J$
induces a metric on
$\De^J$
and on every subcomplex
$K\sub\De^J$.
If
$J$
is finite then
$|J|-1=\dim\De^J$
is the (combinatorial) dimension of
$\De^J$.

For every simplicial polyhedron K, there is the canonical
embedding
$u:K\to\De^J$,
where
$J$
is the vertex set of
$K$,
which is affine on every simplex. Its image
$K'=u(K)$
is called the {\em uniformization} of
$K$.
The (combinatorial) dimension of
$K$
is the maximal dimension of its simplices.

Given a vertex
$v\in K$,
its {\em star}
$\ov\st_v\sub K$
consists of all simplices of
$K$
containing
$v$.
The {\em open star}
$\st_v$
of
$v$
is the star without faces which miss
$v$.
If a simplicial polyhedron
$K$
is uniform, then the open star
$\st_v$
of any vertex
$v\in K$
is an open neighborhood of
$v$.

There are several ways to triangulate the product of
simplicial complexes. One needs for that to choose some
ordering of simplices. Since the barycentric triangulation
is canonically ordered, we prefer to use the following
construction which can be found e.g. in the forthcoming book
\cite{BS3}. Given an index set
$J$,
we denote by
$\ov{\ba}\De^J$
the barycentric subdivision of
$\De^J$,
that is a simplicial complex isometric to
$\De^J$.
The vertices of this complex are barycenters
of all simplices (including 0-dimensional).
The simplices of this complex are convex hulls of all
sets
$S$
of vertices with the property: if
$s$, $s'\in S$
then they are barycenters of two simplices,
one of which is contained in the other. If
$J$
is finite then the covering of
$\De^J$
by the open stars of the vertices of its barycentric
subdivision is
$|J|$-colored:
as the color of a star
$\st_v$
we can take the dimension of the face for which
$v$
is the barycenter.

Now, we describe the barycentric triangulation of the product
of two simplices
$\De^{J_1}$, $\De^{J_2}$.
Regard this product
$\De^{J_1}\times\De^{J_2}$
as a complex with faces that are products of standard simplices.
Let
$b_1$, $b_2$
be barycenters of simplices
$S_1\sub\De^{J_1}$, $S_2 \sub\De^{J_2}$.
We call the point
$(b_1,b_2)\in\De^{J_1}\times\De^{J_2}$
the barycenter of the face
$S_1\times S_2.$
We define the {\em barycentric subdivision} of
$\De^{J_1}\times\De^{J_2}$
as a simplicial complex
$\De^{J_1}\ov{\times}_s\De^{J_2}$,
isometric to
$\De^{J_1}\times\De^{J_2}$,
in the following way. The vertices of
this complex are the barycenters of all faces
(including 0-dimensional). The simplices of this complex
are convex hulls of all sets
$S$
of vertices with the property: if
$s$, $s'\in S$
then they are barycenters of two faces, one of which is contained in the
other.

We denote the uniformization of
$\De^{J_1}\ov{\times}_s\De^{J_2}$
by
$\De^{J_1}\times_s\De^{J_2}$
and call it the {\em barycentric triangulation} of
the product
$\De^{J_1}\times\De^{J_2}$.

The canonical bijection
$\phi:\De^{J_1}\times\De^{J_2}\to\De^{J_1}\times_s\De^{J_2}$
is called the {\em barycentric triangulation map.}

Now, let $K_1$, $K_2$ be uniform simplicial polyhedra,
which we can identify with subpolyhedra of the standard simplices
$\De^{J_1}$, $\De^{J_2}$;
we define the {\em barycentric triangulation} of
$K_1\times K_2$
as $\phi(K_1\times K_2)$.

Note that the covering of
$K_1\times_sK_2$
by open stars of its vertices is
$(n_1+n_2+1)$-colored,
where
$n_1=\dim K_1$, $n_2=\dim K_2$:
one can take as the color of a star
$\st_v$
the dimension of the minimal simplex of
$\De^{J_1\cup J_2}$
containing
$v$.

We give without proof the following two technical lemmas,
(for the  proof see \cite{BS3}).

\begin{lem}\label{lem:bilip} Let
$K_1$, $K_2$
be uniform simplicial polyhedra having both
finite dimensions. The barycentric triangulation map
$$\phi:K_1\times K_2\to K_1\times_sK_2$$
is bilipschitz with bilipschitz constant depending only on dimensions of
$K_1$, $K_2$.
\qed
\end{lem}

\begin{lem}\label{lem:star} Let
$K_1$, $K_2$
be uniform simplicial polyhedra,
$K=K_1\times_sK_2$.
For every vertex
$v\in K$,
there are vertices
$v_1\in K_1$, $v_2\in K_2$
such that
$\phi^{-1}(st_v)\subset st_{v_1}\times st_{v_2}$.
\qed
\end{lem}

Recall that the {\em nerve}
$\cN(\cU)$
of a covering
$\cU=\{ U_j\}_{j\in J}$
is the simplicial complex whose vertices are the members
of the covering and a collection of vertices spans a
simplex if and only if the corresponding members of
$\cU$
have nonempty intersection. We assume that the nerves we
consider are uniform unless the opposite
is explicitly stated; moreover, we can regard
$\cN(\cU)$
as subcomplex of
$\De^J$.

Let
$\cU$
be an open locally finite covering of a metric space
$Z$,
so that no member of it coincides with
$Z$.
The barycentric map
$p:Z\to\cN(\cU)$
is defined as follows. For every
$j\in J$
we put
$q_j(z)=\dist(z,Z\sm U_j)$.
Since the covering is open,
$\sum_{j\in J} q_j(z)>0$.
Since no member of
$\cU$
coincides with
$Z$
and the covering is locally finite,
$\sum_{j\in J} q_j(z)<\infty$
for every point
$z\in Z$.
Then we define the coordinate functions of the map
$p:Z\to\R^J$
by
$p_j(z)=q_j(z)/\sum_{j\in J} q_j(z)$.
Clearly
$p(Z)\sub\cN(\cU)$.
For every vertex
$v\in\cN(\cU)$
the preimage
$p^{-1}(\st_v)$
coincides with corresponding element of the covering.

Suppose in addition that the multiplicity
$m(\cU)\le m+1$
is finite and
$L(\cU)\ge d>0$.
Then, it is not difficult to prove (see e.g. \cite{BS1})
that the map
$p$
is Lipschitz with  Lipschitz constant
$$\Lip(p)\le(m+2)^2/d.$$

\begin{lem}\label{lem:prod}
Let $X_1$, $X_2$ be metric spaces.
 Let
$\cU_1$
be an open covering of
$X_1$, $\cU_2$
an open covering of
$X_2$,
both with multiplicities at most
$n_1+1$
and
$n_2+1$
respectively. Then there exists an
$(n_1+n_2+1)$-colored
open covering
$\cW$
of
$X_1\times X_2$
with
$$L^\times(\cW)\ge q(L(\cU_1), L(\cU_2))\quad\text{and}
\quad\mesh^\times(\cW)\le (\mesh(\cU_1),\mesh(\cU_2)),$$
where the constant
$q$
depends only on
$n_1$, $n_2$.
\end{lem}

\begin{proof} According to the definitions of
$L^\times$
and
$\mesh^\times$,
rescaling the metric of
$X_1\times X_2$
in one of the factors does not change the problem, i.e.
if we construct a required covering for say
$aX_1\times X_2$,
then the rescaling the metric back gives a required covering for
$X_1\times X_2$.

So, by rescaling of the metric, we can assume that
$L(\cU_2)=\frac{(n_2+2)^2}{(n_1+2)^2}L(\cU_1)$.
Denote
$$\la:=\frac{(n_1+2)^2}{L(\cU_1)}=\frac{(n_2+2)^2}{L(U_2)}.$$
Let
$K_i$
be the nerve of
$U_i$
and
$p_i:X_i\to K_i$
the corresponding barycentric map,
$i=1,2$.
So we have
$$\Lip(p_i)\le\la.$$

Then for the map
$p:X_1\times X_2\to K_1\times K_2$,
defined by
$p(x_1, x_2)=(p_1(x_1), p_2(x_2))$,
we have
$\Lip(p)\le\la$.

By Lemma~\ref{lem:bilip}, the barycentric triangulation map
$$\phi:K_1\times K_2\to K_1\times_sK_2$$
is Lipschitz with Lipschitz
constant
$c$,
depending only on $n_1$, $n_2$.
Consider the covering
$\cV$
of
$K_1\times_sK_2$
by open stars of its vertices that is
$(n_1+n_2+1)$-colored.
Since the polyhedron
$K_1\times_sK_2$
is uniform, the Lebesgue number of this covering
is bounded below by a positive constant
$l$
depending only on
$\dim (K_1\times_sK_2)$.

Now let
$\cW$
be the open covering of
$X_1\times X_2$
by preimages of elements of
$\cV$
under the map
$\phi\circ p$.
Then
$\cW$
is
$(n_1+n_2+1)$-colored
and
$$L(\cW)\ge\la^{-1}c^{-1}l=\frac{l}{c(n_1+2)^2}L(\cU_1)
   =\frac{l}{c(n_2+2)^2}L(\cU_2),$$
in particular,
$L^\times(\cW)\ge q(L(\cU_1),L(\cU_2))$
with
$q=\frac{l}{c\max\{(n_1+2)^2,(n_2+2)^2\}}$.

By Lemma~\ref{lem:star} for every
$W\in\cW$
there are
$U_1\in\cU_1$, $U_2\in\cU_2$
so that
$W\sub U_1\times U_2$,
thus
$\mesh^\times(\cW)\le(\mesh\cU_1,\mesh\cU_2)$.
The lemma follows.
\end{proof}

We denote the covering from this lemma by
$\cW=\cU*\cV$.
In the following corollary, we use the natural extension of our
notations
$\mesh^\times$
and
$L^\times$
to the case of more than two factors.

\begin{cor}\label{cor:prod} Let
$\cU_i$
be an open covering with multiplicity at most
$n_i+1$
of a metric space
$X_i$, $i\in\{1,\dots,k\}$.
We let
$\cW_i=\cW_{i-1}*\cU_i$
be the covering of
$X_1\times\dots\times X_i$
for
$i=1,\dots,k$,
where
$\cW_1=\cU_1$.
Then
$\cW=\cW_k$
is an open
$m$-colored
covering of
$X_1\times\dots\times X_k$
with
$m=n_1+\cdots +n_k+1$,
$$L^\times(\cW)\ge q(L(\cU_1),\dots,L(\cU_k))$$
and
$$\mesh^\times(\cW)\le(\mesh(\cU_1),\dots,\mesh(\cU_k)),$$
where the constant
$q$
depends only on
$n_1,\dots,n_k$.
\end{cor}

\begin{proof} Arguing as at the beginning of the proof
above and rescaling of the metrics of the factors,
we can assume that
$L(\cU_1)=\dots=L(\cU_k)=r$.

Note that the property
$L^\times(\cU)\ge(r,\dots ,r)$
for a covering
$\cU$
of a product and the property
$L(\cU)\ge r$
are equivalent. Using this and proceeding by induction,
we obtain
$L^\times(\cW)\ge q(r,\dots,r)$
for some constant
$q$
depending only on
$n_1,\dots,n_k$.

By Lemma~\ref{lem:star} for every
$W\in\cW_i$
there are
$W'\in\cW_{i-1}$, $U\in\cU_i$
so that
$W\sub W'\times U$
for every
$i=2,\dots, k$.
Hence
$$\mesh^\times(\cW)\le(\mesh\cU_1,\dots,\mesh\cU_k).$$
\end{proof}

We shall use the following

\begin{lem}\label{lem:Rcov}
Let
$X$, $X_1$, $X_2$
be metric spaces.

\begin{itemize}
\item[(1)] Assume that
$\ldim(X\times[0,1])\le k$.
Then there exist $\tau_0>0$, $\si>1$
so that for every
$0<\tau<\tau_0$
and
$\tau'>0$
there exists a
$(k+1)$-colored
open covering
$\cU$
of
$X\times\R$
with
$\mesh^\times(\cU)\le(\si\tau,\si\tau')$
and
$L^\times(\cU)\ge(\tau,\tau')$;

\item[(2)] Assume that the spaces of the family
$aX_1\times bX_2$
have
$\ldim\le k$
uniformly in
$a$, $b\ge 1$.
Then there exist
$\tau_0>0$, $\si>1$
so that for every
$0<\tau_1,\tau_2<\tau_0$
and
$\tau'_1$, $\tau'_2>0$
there exists a
$(k+3)$-colored
open covering
$\cW$
of
$Z=X_1\times X_2\times\R\times\R$
with
$\mesh^\times(\cW)\le\si\tau$,
$L^\times(\cW)\ge\tau$,
where
$\tau=(\tau_1,\tau_2,\tau_1',\tau_2')$.
\end{itemize}
\end{lem}

\begin{proof} (1) First, we show that
$\ldim(X\times\R)\le k$.
To this end, we represent
$X\times\R=Z_1\cup Z_2$,
where
$Z_1=X\times\cup_{k\in\Z}(k,k+3/4)$,
$Z_2=X\times\cup_{k\in\Z}(k-1/2,k+1/4)$.
Taking two appropriate qualified coverings of these
sets and applying Proposition~\ref{prop:glue}, we see that
$\ldim(X\times\R)\le k$.
Then taking an appropriate covering of
$X\times\R$
and applying a homothety in the
$\R$
factor, we find a desired covering.

(2) Consider the 2-colored covering
$\cU=\cU^1\cup\cU^2$
of
$\R$,
where
$\cU^1=\set{(2k, 2k+2)}{$k\in\Z$}$,
$\cU^2=\set{(2k+1, 2k+3)}{$k\in\Z$}$.
Clearly,
$L(\cU)=1/2$, $\mesh(\cU)=2$.

There is a constant
$\de\in(0,1)$
such that for every sufficiently small
$\tau_1$, $\tau_2>0$
there exists a
$(k+1)$-colored
open covering
$\cV$
of
$X\times  Y$
with
$\mesh^\times(\cV)\le(\tau_1,\tau_2)$
and
$L^\times(\cV)\ge\de(\tau_1,\tau_2)$.

For the product covering
$\cW=(\cV*\cU)*\cU$
of
$Z$
we have by Corollary~\ref{cor:prod}
$L^\times(\cW)\ge q(\tau_1,\tau_2,1/2,1/2)$
and
$\mesh^\times(\cW)\le(\tau_1,\tau_2,2,2)$
for some constant
$q>0$
depending only on
$k$, $\de$.
Then applying homotheties independently in two
$\R$-factors
and changing the constants, we obtain a desired covering.
\end{proof}

\subsection{Locally self-similar spaces}

Let
$\la\ge 1$ and
$R>0$
be given. A map
$f:Z\to Z'$
between metric spaces is
$\la$-{\em quasi-homothetic
with coefficient}
$R$
if for all
$z$, $z'\in Z$,
we have
$$R|zz'|/\la\le|f(z)f(z')|\le\la R|zz'|.$$

A metric space
$Z$
is {\em locally self-similar},
if there is
$\la\ge 1$
such that for every sufficiently large
$R>1$
and every
$A\sub Z$
with
$\diam A\le 1/R$,
there is a
$\la$-quasi-homothetic
map
$f:A\to Z$
with coefficient
$R$.

It is proved in \cite{BL} that the linearly controlled
dimension
$\ldim Z$
of every compact locally self-similar
metric space
$Z$
is finite and coincides with
$\dim Z$, $\ldim Z=\dim Z$.

We shall use the following facts obviously implied
by the definition of a quasi-homothetic map.

\begin{lem}\label{lem:distort} Let
$h:Z\to Z'$
be a
$\la$-quasi-homothetic
map with coefficient
$R$.
Let
$V\sub Z$
and let
$\wt\cU$
be an open covering of
$h(V)$, $\cU=h^{-1}(\wt\cU)$.
Then
\begin{itemize}
\item[(1)] $R\mesh(\cU)/\la\le\mesh(\wt\cU)
           \le\la R\mesh(\cU)$;

\item[(2)] $\la R\cdot L(\cU)\ge L(\wt\cU)
    \ge R\cdot L(\cU)/\la$,
where
$L(\cU)$
is the Lebesgue number of
$\cU$
as a covering of
$V$.\qed
\end{itemize}
\end{lem}

\begin{pro}\label{prop:ssim} If compact metric spaces
$X$, $Y$
are locally self-similar
then the spaces of the family
$a X\times bY$
have
$\ldim\le n$, $n=\dim(X\times Y)$,
uniformly in
$a$, $b\ge 1$.
\end{pro}

\begin{proof} By the remark above,
$\ldim X=N$
is finite. Furthermore, the space
$Z=X\times Y$
is also compact and locally self-similar. Thus
$\ldim Z=n$
is finite.

It suffices to prove that there is a constant
$\de\in(0,1)$
such that for every sufficiently small
$\tau>0$
and for every
$\al$, $\be\in(0,1]$
there exist
$(n+1)$-colored
open coverings
$\cU$, $\cU'$
of
$Z$
with
$\mesh^\times(\cU)\le(\al\tau,\tau)$,
$L^\times(\cU)\ge\de(\al\tau,\tau)$;
$\mesh^\times(\cU')\le (\tau,\be\tau)$
and
$L^\times(\cU')\ge\de(\tau,\be\tau)$.

Let
$\de_1$, $\de_2$, $\de_3$
be the constants from the definition of
$\ldim$
for
$X$, $Y$, $Z$
respectively, and let
$\de'=\min\{\de_1,\de_2,\de_3\}$.
We can assume that the constant
$\la$
from the definition of self-similarity
is common for both
$X$
and
$Y$.

We explain how to construct the covering
$\cU$,
the covering
$\cU'$
is constructed similarly. Fix a positive
$\tau<\min\{\de'/\la,\diam Y\}$.
For every
$\al\in(0,1]$,
we construct a covering with
$\mesh^\times\le(\alpha\tau,\tau)$
and
$L^\times\ge\de(\alpha\tau,\tau)$,
where
$\de=(\de'/4\la^2)^{N+1}/2$.

The argument is similar to that in \cite[Theorem~1.1]{BL}.
We fix an
$(N+1)$-colored
open covering
$\cV''$
of
$X$
with
$\mesh(\cV'')\le\al/\la$,
$L(\cV'')\ge\de'\al/\la$
and put
$\cV=\set{V\times Y}{$V\in \cV''$}$.
Then
$\cV$
is an
$(N+1)$-colored
open covering of
$Z$
with
$L(\cV)\ge\de'\al/\la$,
and we assume that
$\cV$
is colored by the set $S=\{0,\dots,N\}$,
$\cV=\cup_{a\in S}\cV^a$.

For every
$V\in\cV$,
consider the slightly smaller subset
$V'=B_{-\de'\al/2\la}(V)$.
Then, the sets
$Z_a=\cup_{V\in\cV^a}V'\sub Z$, $a\in S$,
cover
$Z$, $Z=\cup_{a\in S}Z_a$,
because
$L(\cV)\ge\de'\al/\la$.
The idea is to construct a family of open
$(n+1)$-colored
coverings
$\cU_a$
of
$Z_a$, $a\in S$,
which is separated with and qualified by the scale function
$i:S\to\{0,\dots,N\}$, $i(a)=N-a$,
and then to construct the desired
$\cU$
using Proposition~\ref{prop:glue}. The construction of
$\cU_a$
is based on self-similarity of
$X$.

Using that
$\ldim(X\times Y)=n$
and assuming that
$\tau$
is sufficiently small, we find for every
$a\in S$
an
$(n+1)$-colored
covering
$\wt\cU_a$
of
$X\times Y$
with
$\mesh(\wt\cU_a)\le(\de'/4\la^2)^a\tau/4$
and
$L(\wt\cU_a)\ge\de'(\de'/4\la^2)^a\tau/4$.

Given
$V\in\cV$,
we fix a map
$h_V=(h_V^1,\id):V\to X\times Y$,
where
$h_V^1$
is
$\la$-quasi-homothetic
with coefficient
$R=\la/\al$,
and put
$\wt V=h_V(V')$.
Now, for every
$a\in S$, $V\in\cV^a$
consider the family
$\wt\cU_{a,V}=\set{\wt U\in\wt\cU_a}{$\wt V\cap\wt U\neq\es$}$
which is obviously an
$(n+1)$-colored
covering of
$\wt V$.
Then,
$$\cU_{a,V}=\set{h_V^{-1}(\wt U)}{$\wt U\in\wt\cU_{a,V}$}$$
is an open
$(n+1)$-colored
covering of
$V'$.

Note that
$U=h_V^{-1}(\wt U)$
is contained in
$V$
for every
$\wt U\in\wt\cU_{a,V}$
because
$\mesh^\times\wt U\le(\tau/4,\tau/4)$,
therefore
$\mesh^\times U\le(\al\tau/4,\tau/4)$
by Lemma~\ref{lem:distort}, and hence
$U\sub B_{(\al\tau/2,\diam Y)}(V')\sub V$
by the choice of
$\tau$.
Thus the family
$\cU_{a,V}$
is contained in
$V$.
Now, the family
$\cU_a=\cup_{V\in\cV^a}\cU_{a,V}$
covers the set
$Z_a$
of the color
$a$,
and it has the following properties

\begin{itemize}
\item[(1)] for every
$a\in S$,
the covering
$\cU_a$
is
$(n+1)$-colored;

\item[(2)] $L(\cU_a)^\boxtimes\ge
    4\mesh^\boxtimes\cU_{a+1}$
for every
$a\in S$, $a\le N-1$;

\item[(3)] $\mesh^\times(\cup_{a\in S}\cU_a)
   \le(\alpha\tau,\tau)$,
  $L^\times(\cup_{a\in S}\cU_a)
  \ge(\de'/4\la^2)^{N+1}(\al\tau,\tau)$.
\end{itemize}

Indeed, distinct
$V_1$, $V_2\in\cV^a$
are disjoint and thus any
$U_1\in\cU_{a,V_1}$, $U_2\in\cU_{a,V_2}$
are disjoint because
$U_1\sub V_1$, $U_2\sub V_2$.
This proves (1).

Applying Lemma~\ref{lem:distort}, we see that
$\mesh^\times(\cU_a)\le\mesh(\wt\cU_a)(\alpha,1)$
and
$L^\times(\cU_a)\ge L(\wt\cU_a)(\al/\la^2,1/\la^2)$
for every
$a\in S$.
These estimates together with the estimates on
$\mesh(\wt\cU_a)$, $L(\wt\cU_a)$
yield (3), and together with the inequalities
$$4\la^2\mesh(\wt\cU_{a+1})\le\de'(\de'/4\la^2)^a\tau/4
   \le L(\wt\cU_a)$$
for every
$a\in S$, $a\le N-1$,
yield (2).

In view (1) and (2), the family of coverings
$\cU_a$, $a\in S$,
satisfies the condition of Proposition~\ref{prop:glue}.
Applying this proposition and using (3), we obtain
an open
$(n+1)$-colored covering
$\cU$
of
$X\times Y$
with
$\mesh^\times(\cU)\le(\al\tau,\tau)$
and
$L^\times(\cU)\ge\de(\al\tau,\tau)$,
where
$\de=(\de'/4\la^2)^{N+1}/2$.
\end{proof}

\subsection{`Hyperbolic cone' map}
Let
$Z$
be a bounded metric space. Assuming that
$\diam Z>0$,
we put
$\mu=\pi/\diam Z$
and note that
$\mu|zz'|\in[0,\pi]$
for every
$z$, $z'\in Z$.
Recall that the hyperbolic cone
$\cone(Z)$
over
$Z$
is the space
$Z\times[0,\infty)/Z\times\{0\}$
with metric defined as follows. Given
$x=(z,t)$, $x'=(z',t')\in\cone(Z)$
we consider a triangle
$\ov o\,\ov x\,\ov x'\sub\hyp^2$
with
$|\ov o\,\ov x|=t$, $|\ov o\,\ov x'|=t'$
and the angle
$\angle_{\ov o}(\ov x,\ov x')=\mu|zz'|$.
Now, we put
$|xx'|:=|\ov x\,\ov x'|$.
In the degenerate case
$Z=\{\pt\}$,
we define
$\cone(Z)=\{\pt\}\times[0,\infty)$
as the metric
product. The point
$o=Z\times\{0\}\in\cone(Z)$
is called the {\em vertex} of
$\cone(Z)$.

We let
$h:Z\times[0,\infty)\to\cone(Z)$
be the canonical projection, and
$\de=\de_{\hyp^2}$
the hyperbolicity constant of
$\hyp^2$.

It is well known that the hyperbolic cone
$\cone(Z)$
is a hyperbolic space which satisfies the
$\de$-inequality
w.r.t. the vertex
$o$.
Furthermore, there is a canonical inclusion
$Z\sub\di\cone(Z)$,
and the metric of
$Z$
is visual w.r.t. the base point
$o$
and the parameter
$e$,
$$e^{-(\xi|\xi')_o-c_0}\le d(\xi,\xi')\le e^{-(\xi|\xi')_o+c_0}$$
for some
$c_0\ge 0$
and all
$\xi$, $\xi'\in Z$.
In general,
$\cone(Z)$
is not geodesic, however, for every point
$z\in\cone(Z)$
there is a geodesic segment
$oz\sub\cone(Z)$
and every
$\xi\in Z$
is represented by a geodesic ray
$o\xi\sub\cone(Z)$.

The following lemma is similar to \cite[Lemma~5.1]{BoS}.

\begin{lem}\label{lem:grprest} Let
$X$
be a
$\de$-hyperbolic space. Given
$o$, $z_1$, $z_2\in X$
and
$x_1\in oz_1$, $x_2\in oz_2$,
we have
$$\min\{|ox_1|,|ox_2|,(z_1|z_2)_o\}-2\de\le(x_1|x_2)_o
\le\min\{|ox_1|,|ox_2|,(z_1|z_2)_o\}+2\de.$$
\end{lem}

\begin{proof} Applying the
$\de$-inequality
twice, we obtain
\begin{eqnarray*}
(z_1|z_2)_o\ge\min\{(z_1|x_1)_o, (x_1|x_2)_o, (x_2|z_2)_o\}-2\de \\
 =\min\{|ox_1|,(x_1|x_2)_o,|ox_2|\}-2\de=(x_1|x_2)_o-2\de.
\end{eqnarray*}
Similarly we have
$(x_1|x_2)_o\ge\min{|ox_1|,(z_1|z_2)_o,|ox_2|}-2\de$.
The lemma follows.
\end{proof}

\begin{cor}\label{cor:GromProd} Let
$X$
be a
$\de$-hyperbolic space.
Suppose that points
$x_1\in o\xi_1$, $x_2\in o\xi_2$
satisfy
$|x_1o|=|x_2o|$
and
$|x_1x_2|>4\de$,
where
$\xi_1$, $\xi_2\in\di X$.
Then
$$(x_1|x_2)_o-2\de\le(\xi_1|\xi_2)_o\le(x_1|x_2)_o+2\de.$$
\end{cor}

\begin{proof} In this case
$(x_1|x_2)_o<|ox_1|-2\de=|ox_2|-2\de$.
Hence by the previous lemma,
$\min\{|ox_1|,|ox_2|,(z_1|z_2)_o\}=(z_1|z_2)_o$
for each
$z_1\in(x_1,\xi_1)$, $z_2\in(x_2,\xi_2)$.
Now, the statement follows from that lemma.
\end{proof}

We put
$c=c_0+2\de$.

\begin{claim}\label{cl:h1} For every
$R>0$,
the map
$h|_{Z\times\{R\}}:Z\times\{R\}\to\cone(Z)$
possesses the following properties

\begin{itemize}
\item[(1)] $|h(z,R)h(z',R)|\ge 2D$
for every
$D>0$
and each
$z$, $z'\in Z$
with
$|zz'|\ge e^{-R+D+c}$;

\item[(2)] $|h(z,R)h(z'R)|\le 2D$
for every
$D>2\de$
and each
$z$, $z'\in Z$
with
$|zz'|\le e^{-R+D-c}$.
\end{itemize}
\end{claim}

\begin{proof} (1) We identify
$z=\xi$, $z'=\xi'\in\di\cone(Z)$
and denote
$h(z,R)=x$, $h(z',R)=x'$.
Thus
$$e^{-R+D+c}\le|\xi\xi'|\le e^{-(\xi|\xi')_o+c_0},$$
and we obtain
$(\xi|\xi')_o\le R-D-2\de$.
By Lemma~\ref{lem:grprest},
$$(x|x')_o\le\min\{R,(\xi|\xi)_o\}+2\de\le R-D.$$
Hence,
$|xx'|=2(R-(x|x')_o)\ge 2D$.

(2) Using the notations as above, we obtain
$e^{-(\xi|\xi')_o-c_0}\le|\xi\xi'|\le e^{-R+D-c}$.
Hence,
$(\xi|\xi')_o\ge R-D+2\de$.
In the case
$|xx'|\le 4\de$,
we see that
$|xx'|<2D$
by the assumption
$D>2\de$.
In the case
$|xx'|>4\de$,
we apply Corollary~\ref{cor:GromProd} and obtain
$(x|x')_o\ge(\xi|\xi')_o-2\de\ge R-D$.
Thus
$|xx'|=2(R-(x|x')_o)\le 2D$.
\end{proof}

The proof of the following claim is similar, and we leave it
as an exercise to the reader. We only mention that at some
point one should use monotonicity of the Gromov product,
$(x|x')_o\le(\xi|\xi')_o$
for every
$x\in o\xi$, $x'\in o\xi'$.

\begin{claim}\label{cl:h2} For every
$R>0$,
the map
$h|_{Z\times\{R\}}:Z\times\{R\}\to\cone(Z)$
possesses the following properties

\begin{itemize}
\item[(1)] $|zz'|\ge e^{-R+D-c}$
for every
$D>2\de$
and each
$x=h(z,R)$, $x'=h(z',R)$
with
$|xx'|\ge 2D$;

\item[(2)] $|zz'|\le e^{-R+D+c}$
for every
$D>0$
and each
$x=h(z,R)$, $x'=h(z',R)$
with
$|xx'|\le 2D$.
\end{itemize}
\qed
\end{claim}

\begin{claim}\label{cl:h3} (1) Given
$l>0$, $R\ge R_1\ge l$,
we let
$r=e^{-R_1+3l/2+c}$.
Then
$$h(B_{r}(x)\times[R-l, R+l])\sups B_{l}(h(x,R));$$

(2) given $l>2\de$, $t>0$, $R_2\ge R\ge t$,
we let
$r=e^{-R_2+l-c}$.
Then
$$h(B_{r}(x)\times [R-t, R+t])\sub B_{3t+2l}(h(x,R)).$$
\end{claim}

\begin{proof} (1) Given
$y'\in \cone(Z)$
with
$|y'h(x,R)|<l$,
we represent
$y'=h(y,R_0)$.
Using the definition of distances in
$\cone(Z)$,
we obtain
$R_0\in[R-l,R+l]$
and
$|h(y,R')h(x,R')|<l$,
where
$R'=\min\{R_0,R\}$.
Then by Claim~\ref{cl:h2}(2), we have
$$|yx|\le e^{-R'+l/2+c}\le e^{-R+l+l/2+c}\le e^{-R_1+3l/2+c}.$$
Claim~\ref{cl:h3}(1) follows.

(2) Assume that
$(y,R_0)\in B_{r}(x)\times[R-t,R+t]$.
Then
$R_0\in[R-t,R+t]$
and
$|yx|<e^{-R_2+l-c}\le e^{-R_0+t+l-c}$.
By Claim~\ref{cl:h1}(2), we have
$|h(y,R_0)h(x,R_0)|<2t+2l$.
Then by the triangle inequality,
$$|h(x,R)h(y,R_0)|<3t+2l.$$
Claim~\ref{cl:h3}(2) follows.
\end{proof}

\begin{claim}\label{cl:h4}
(1) Given
$R_1\ge R>0$, $l>8\de$,
we let
$r_1=e^{-R_1-l/4-c}$.
Then
$h^{-1}(B_{l}(h(x,R)))\sups
    B_{r_1}(x)\times[R-l/2,R+l/2]$;

(2) given $l>0$, $R\ge R_2\ge l$,
we let
   $r_2=e^{-R_2+3l/2+c}$.
Then
 $h^{-1}(B_l(h(x,R)))\sub
    B_{r_2}(x)\times[R-l,R+l]$.
\end{claim}

\begin{proof} (1) Assume that
$(y,R_0)\in B_{r_1}(x)\times[R-l/2,R+l/2]$.
Then
$R_0\in [R-l/2,R+l/2]$
and
$|xy|<e^{-R_1-l/4-c}\le e^{-R_0+l/4-c}$.
Using Claim~\ref{cl:h1}(2), we obtain
$|h(y,R_0)h(x,R_0)|<l/2$.
Then by the triangle inequality,
$|h(y,R_0)h(x,R)|<l$.
Claim~\ref{cl:h4}(1) follows.

(2) Claim~\ref{cl:h4}(2) follows from Claim~\ref{cl:h3}(1).
\end{proof}

\section{Proofs}

\begin{pro}\label{prop:asprod<} Let
$Z_1$, $Z_2$
be bounded metric spaces such that the spaces
of the family
$\set{aZ_1\times bZ_2}{$a\ge1,b\ge1$}$
have
$\ldim\le k$
uniformly. Then
$$\asdim(\cone(Z_1)\times\cone(Z_2))\le k+2.$$
\end{pro}

We need some preparation for the proof. We use
the following notations

\begin{eqnarray*}
\lo {P_0}(T)&=&Z_1\times[0,T]\times Z_2\times [0,T],\\
\lo {P_1}(T)&=&Z_1\times[T,\infty]\times Z_2\times [0,T],\\
\lo {P_2}(T)&=&Z_1\times[0,T]\times Z_2\times [T,\infty],\\
\lo {P_3}(T)&=&Z_1\times[T,\infty]\times Z_2\times [T,\infty],
\end{eqnarray*}
and
$P_i=(h_1,h_2)(\lo{P_i})$,
where
$h_i:Z_i\times[0,\infty)\to\cone Z_i$
is the canonical projection,
$i=1,2$.

\begin{lem}\label{lem:lmk} For every
$L>0$
there exist
$T_3$, $M>0$
so that for every
$T>T_3$
there exists a
$(L,M,k+3)$-covering of
$P_3(T)$.
\end{lem}

\begin{proof} Let
$\si>1$
be the constant from Lemma~\ref{lem:Rcov}(2).
We fix a constant
$L>0$,
put
$H=512\si^4L$
and for every integer
$m$, $n\ge 0$
consider the product

$$A_{m,n}=Z_1\times Z_2\times[Hm,Hm+H]\times[Hn,Hn+H].$$

For some
$T>0$,
we first construct a
$(k+3)$-colored
covering
$\cU$
of
$\lo {P_3}(T)$,
so that for all sufficiently large integers
$m$, $n$
the following holds
$$L^\times(\cU)|_{A_{m,n}}
  \ge(e^{-Hm+3 L/2+c},e^{-Hn+3 L/2+c},L,L);$$
$$\mesh^\times(\cU)|_{A_{m,n}}\le(e^De^{-Hm},e^De^{-Hn},H/4,H/4)$$
for some
$D>0$,
where the constant
$c>0$
is defined just before Claim~\ref{cl:h1}.

It is convenient to use the following notations
$c_1=2e^{3L/2+c}$, $b=e^{-H}$.
Using Lemma~\ref{lem:Rcov}(2), we fix a sufficiently large
$N$
so that for each integers
$m$, $n\ge N$
and for every
$i=0,1,2,3$
there exists a
$(k+3)$-colored
covering
$\wt\cU^i_{m,n}$
of
$Z_1\times Z_2\times\R\times\R$
with
$L^\times(\wt\cU^i_{m,n})\ge(4\si )^i(c_1b^{m-i},c_1b^{n-i},2L,2L)$
and
$\mesh^\times(\wt\cU^i_{m,n})\le\si(4\si)^i(c_1b^{m-i},c_1b^{n-i},2L,2L)$.
For
$m$, $n\ge N$,
we consider the subfamily
$\cU^i_{m,n}$
of
$\wt\cU^i_{m,n}$
consisting of all members that intersect
$A_{m,n}$.

Next, we define the function
$i:\N\times\N\to\{0,1,2,3\}$
by
$i(m,n):=m\mod2+2(n\mod2)$
and put
$\cU_{m,n}:=\cU_{m,n}^{i(m,n)}$.
Now, we show that the family of coverings
$\cU_{m,n}$
of
$A_{m,n}$, $m$, $n\ge N$,
is separated with and qualified by the scale function
$i$.

We have
$\si(4\si)^32L=H/4<H/2$
by the choice of the constants. It follows from the
estimate on
$\mesh^\times(\wt\cU_{m,n}^i)$
that any member
$U\in\cU_{m,n}$
is disjoint with any
$U'\in\cU_{m',n'}$,
if
$i(m,n)=i(m',n')$
and
$(m,n)\neq(m',n')$.
That is, the family
$\cU_{m,n}$, $m$, $n\ge N$,
is separated with the scale function
$i$.
In particular, if
$U\cap U'\neq\es$
for some
$U\in\cU_{m,n}$, $U'\in\cU_{m',n'}$
then the pairs
$(m,n)$
and
$(m',n')$
are adjacent, i.e.,
$|m'-m|\le 1$
and
$|n'-n|\le 1$.

Assume that pairs
$(m',n')$, $(m,n)$
are adjacent and
$i(m',n')>i(m,n)$.
We have
$(4\si)^{i+1}2L=4\si(4\si)^i2L$
and
$(4\si)^{i+1}c_1b^{m-(i+1)}>4\si(4\si)^ic_1b^{m-i}$,
since
$b<1$.
These inequalities together with the estimates on
$L^\times(\wt\cU_{m,n}^i)$, $\mesh^\times(\wt\cU_{m,n}^i)$,
show that
$$L^\boxtimes(\cU_{m',n'})\ge 4\mesh^\boxtimes(\cU_{m,n}).$$
It follows that these coverings are qualified by the scale function
$i$.
So applying Proposition~\ref{prop:glue} to this family of coverings,
we obtain a
$(k+3)$-colored
covering
$\cU$
of
$\ov P_3(T)$, $T=HN$,
for which the following holds

\begin{eqnarray*}
 L^\times(\cU)|_{A_{m,n}}
   &\ge&\frac{1}{2}(c_1e^{-Hm},c_1e^{-Hn},2L,2L)\\
   &=&(e^{-Hm+3L/2+c},e^{-Hn+3L/2+c},L,L);\\
  \mesh^\times(\cU)|_{A_{m,n}}
  &\le&\si(4\si)^3
  (c_1e^{-H(m-i(m,n))},c_1e^{-H(n-i(m,n))},2L,2L)\\
  &\le&(e^De^{-Hm},e^De^{-Hn},H/4,H/4)
\end{eqnarray*}
for some
$D>0$.

Now, we estimate the Lebesgue number and the mesh of
the covering
$\cU'=(h_1, h_2)(\cU)$.
Applying Claim~\ref{cl:h3}(1) with
$l=L$, $R_1=H_m$, $R\in[Hm,Hm+H]$
(resp.
$R_1=Hn$, $R\in[Hn,Hn+H]$)
to the map
$h_1$
(resp.
$h_2$),
we obtain
$L(\cU')\ge L$.

Next, we represent the estimate above for
$\mesh^\times(\cU)|_{A_{m,n}}$
as follows
$$\mesh^\times(\cU)|_{A_{m,n}}
  \le(e^{-(Hm+H)+H+D+c-c},e^{-(Hn+H)+H+D+c-c},H/4,H/4).$$
Now, applying Claim~\ref{cl:h3}(2) with
$l=H+D+c$
(note that
$l>2\de=2\de_{\hyp^2}$), $t=H/4$, $R_2=Hm+H$, $R\in[Hm,Hm+H]$
(resp. $R_2=Hn+H$, $R\in[Hn,Hn+H]$)
to the map
$h_1$
(resp.
$h_2$),
we obtain
$\mesh(\cU')\le M=2(3H/4+2(H+D+c))$.

It follows that
$\cU'$
is a
$(L,M,k+3)$-covering
of
$P_3(T)$.
\end{proof}

\begin{lem}\label{lem:lmk12} Given
$L>0$,
there is
$T_2$
so that for every
$T>T_2$
there exist
$(L,M,k+2)$-coverings
of
$P_1(T)$, $P_2(T)$
with
$M>0$
depending on
$T$, $L$.
\end{lem}

\begin{proof} The proof is similar to that of
Lemma~\ref{lem:lmk}. For simplicity, we construct a
desired covering of
$P_1(T)$.
For
$P_2(T)$
the construction is the same with obvious modifications.

We fix
$L>0$,
put
$H=40\si^2L$, $\si$
is the constant from Lemma~\ref{lem:Rcov}(1),
and for every integer
$m\ge 0$
consider the product
$$A_m=Z_1\times  Z_2\times[Hm,Hm+H]\times[0,T].$$
For some
$T>0$,
we first construct a
$(k+2)$-colored
covering
$\cU$
of
$\lo {P_1}(T)$,
so that for all
$m\ge N$
the following holds
$$L(\cU)|_{A_m}\ge(e^{-Hm+3L/2+c},\diam Z_2,L,L),$$
$$\mesh(\cU)|_{A_m}\le(e^De^{-Hm},\diam Z_2,H/4,T+L)$$
for some
$D>0$,
where the constant
$c=c_0+2\de$
is introduced before Claim~\ref{cl:h1}.

For convenience we use notations
$c_1=2e^{3L/2+c}$, $b=e^{-H}$.
Since
$Z_1$
is isometrically embedded into
$Z_1\times Z_2$,
we have
$\ldim Z_1\le\ldim(Z_1\times Z_2)\le k$
and therefore
$\ldim(Z_1\times[0,1])\le k+1$.
Then, using Lemma~\ref{lem:Rcov}(1), we find a sufficiently large
$N$
so that for every integer
$m\ge N$
and every
$i=0,1$
there exists a
$(k+2)$-colored
covering
$\wt\cU^i_{m}$
of
$Z_1\times\R$
with
\begin{eqnarray*}
 L^\times(\wt\cU^i_{m})&\ge&(4\si)^i(c_1b^{m-i},2L)\\
 \mesh^\times(\wt\cU^i_{m})&\le&\si(4\si )^i(c_1b^{m-i},2L).
\end{eqnarray*}
We take its subfamily
$\cU^i_{m}$
consisting of all members that intersect
$\wh A_m=Z_1\times[Hm,Hm+h]$.

We define the function
$i:\N\to\{0,1\}$
by
$i(m):=m\mod2$
and put
$\cU_{m}:=\cU_{m}^{i(m)}$.
Now, we show that the family of the coverings
$\cU_{m}$
of
$\wh A_m$, $m\ge N$,
is separated with and qualified by the scale function
$i$.

We have
$\si(4\si)2L<H/4$
by the choice of the constants. It follows from the
estimate on
$\mesh^\times(\wt\cU_{m}^i)$
that any member
$U\in\cU_{m}$
is disjoint with any
$U'\in\cU_{m'}$,
if
$i(m)=i(m')$
and
$m\neq m'$.
That is, the family
$\cU_{m}$, $m\ge N$,
is separated with the scale function
$i$.
In particular, if
$U\cap U'\neq\es$
for some
$U\in\cU_{m}$, $U'\in\cU_{m'}$
then
$|m'-m|\le 1$.

Assume that
$|m'-m|=1$
and
$i(m')>i(m)$,
i.e.,
$i(m')=1, i(m)=0$.
Using the estimates on
$L^\times(\wt\cU_m^i)$, $\mesh^\times(\wt\cU_m^i)$
and the fact that
$b<1$,
we obtain
$$L^\boxtimes(\cU_{m'})\ge 4\mesh^\boxtimes(\cU_m).$$
Thus these coverings are qualified by
the scale function
$i$.
So, we apply Proposition~\ref{prop:glue} to the family of
the coverings
$\cU_m$, $m\ge N$,
and obtain a
$(k+2)$-colored
covering
$\cU'$
of
$Z_1\times[T,\infty)$, $T=HN$,
with
\begin{eqnarray*}
 L^\times(\cU')|_{A_m}
   &\ge&\frac{1}{2}(c_1e^{-Hm},2L)
   =(e^{-Hm+3L/2+c},L);\\
  \mesh^\times(\cU')|_{A_m}
  &\le& 4\si^2
  (c_1e^{-H(m-i(m))},2L)
  \le(e^De^{-Hm},H/4)
\end{eqnarray*}
for some
$D>0$.

Now, the family
$\cU=\set{U\times Z_2\times [0,T+L]}{$U\in \cU'$}$
is a
$(k+2)$-colored
covering of
$\lo {P_1}(T)$
with required estimates on
$L^\times$-
and
$\mesh^\times$-numbers.

Let us estimate the Lebesgue number and the mesh of
the covering
$\cU''=(h_1, h_2)(\cU)$.
Applying Claim~\ref{cl:h3}(1) with
$l=L$, $R_1=H_m$, $R\in[Hm,Hm+H]$
to the map
$h_1$,
we obtain
$L(\cU'')\ge L$.

Next, we represent the estimate above for
$\mesh^\times(\cU)|_{A_m}$
as follows
$$\mesh^\times(\cU)|_{A_m}
   \le(e^{-(Hm+H)+H+D+c-c},\diam Z_2,H/4,T+L).$$
Now, applying Claim~\ref{cl:h3}(2) with
$l=H+D+c$
(note that
$l>2\de=2\de_{\hyp^2}$), $t=H/4$, $R_2=Hm+2H$, $R\in[Hm,Hm+H]$
to the map
$h_1$
and using that the diameter of
$h_2(Z_2\times[0,T+L])\sub\cone(Z_2)$
is at most
$2(T+L)$,
we obtain
$\mesh(\cU'')\le M=2\max\{3H/4+2(2H+D+c),T+L\}$.

It follows that
$\cU''$
is a
$(L,M,k+2)$-covering
of
$P_1(T)$.
\end{proof}

\begin{lem}\label{lem:lmk0} For every
$T$, $L>0$
there exists a
$(L,M,1)$-covering
of
$P_0(T)$
with
$M>0$
depending only on
$T$
and
$L$.
\end{lem}

\begin{proof} The covering of
$P_0(T)$
consisting of the unique element
$U=h_1(Z_1\times[0,T+2L))\times h_2(Z_2\times[0,T+2L))$
has the required properties.
\end{proof}

\begin{proof}[Proof of Proposition~\ref{prop:asprod<}]
Given
$L_0>0$,
we construct a uniformly bounded covering
$\cU$
of
$\cone(Z_1)\times\cone(Z_2)$
with
$L(\cU)\ge L_0$.

We apply Lemma~\ref{lem:lmk} for
$L=L_0$,
find corresponding constants
$T_3$, $M_3$
and let
$\cU_3$
be a
$(L_0,M_3,k+3)$-covering of
$P_3(T_3)$.
We can assume that
$M_3\ge L_0$.

Next, we apply Lemma~\ref{lem:lmk12} first for
$L=2M_3$
and denote by
$T_2$
the corresponding constant. Thus for
$T_0=\max\{T_2, T_3 \}$
there is a
$(2M_3,M_2,k+2)$-covering
$\cU_2$
of
$P_2(T_0)$,
where the constant
$M_2\ge 2M_3$
depends on
$T_0$, $M_3$.

Then we again apply Lemma~\ref{lem:lmk12} for
$L=2M_2$
and denote by
$T_1$
the corresponding constant. Thus for
$T^*=\max\{T_0, T_1\}$
there is a
$(2M_2,M_1,k+2)$-covering
$\cU_1$
of
$P_1(T^*)$,
where the constant
$M_1\ge 2M_2$
depends on
$T^*$, $M_2$.

Finally, by Lemma~\ref{lem:lmk0}, there is a
$(2M_1,M,1)$-covering
$\cU_0$
of
$P_0(T^*)$,
where
$M>0$
depend on
$M_1$
and
$T^*$.

We have
$$\cone(Z_1)\times\cone(Z_2)
   =P_0(T^*)\cup P_1(T^*)\cup P_2(T_0)\cup P_3(T_3).$$

Now we apply Lemma~\ref{lem:union} consequently to the
sequence of the coverings
$\cU_0$, $\cU_1$, $\cU_2$, $\cU_3$
and obtain a bounded
$(k+3)$-colored
covering
$\cU$
of
$\cone(Z_1)\times\cone(Z_2)$
with
$L(\cU)\ge L_0$.
This shows that
$\asdim(\cone(Z_1)\times\cone(Z_2))\le k+2$.
\end{proof}

\begin{pro}\label{prop:lasdim} Let
$Z$
be a bounded metric space. Then
$$\lasdim\cone(Z)\le\ldim(Z\times[0,1]).$$
\end{pro}

\begin{proof} We can assume that
$\ldim(Z\times[0,1])=k$
is finite. We fix some
$L>4\de$, $\de=\de_{\hyp^2}$.
Let
$\tau_0>0$, $\si>1$
be the constants from Lemma~\ref{lem:Rcov}(1).
We put
$H=40\si^2 L$
and denote
$$A_m=Z\times[Hm,Hm+H].$$
There exists
$N>0$
such that for
$m\ge N$
we have
$e^{-H(m-1)+3L/2+c}<\tau_0$
uniformly in
$L>4\de$,
where as above the constant
$c=c_0+2\de$
is introduced before Claim~\ref{cl:h1}.

Exactly as in the proof of Lemma~\ref{lem:lmk12},
we construct a
$(k+1)$-colored
covering
$\cU$
of
$Z\times[HN,\infty)$
such that for all
$m\ge N$
we have
\begin{eqnarray*}
 L^\times(\cU)|_{A_{m}}&\ge& (e^{-Hm+3L/2+c}, L),\\
 \mesh^\times(\cU)|_{A_{m}}&\le& (e^{D}e^{-Hm},H/4)
\end{eqnarray*}
for some constant
$D>0$
linearly depending on
$L$.

We have
$$\mesh^\times(\cU)|_{A_m}\le (e^{-(Hm+H)+(H+3L/2+D+c)-c},H/4).$$
Applying Claim~\ref{cl:h3} to the canonical map
$h:Z\times[0,\infty)\to\cone(Z)$,
we see that the covering
$h(\cU)=\cU'$
has
$L(\cU')\ge L$
and
$\mesh(\cU')\le(3H/4+2(H+D+3L/2+c))$.

Now for every
$L>4\de$
we have the
$(k+1)$-colored
covering
$\cU'$
of
$h(Z\times[CL,\infty))$,
$L(\cU')\ge L$, $\mesh(\cU')\le CL$
for some constant
$C\ge 1$
independent of
$L$.
Then we cover the ball
$B_{CL}(o)\sub\cone Z$
by the ball
$\cU''=B_{3CL}(o)$
and apply Lemma~\ref{lem:union} to
$\cU'$, $\cU''$
to obtain a covering
$\cW$
of
$\cone Z$
with multiplicity
$m(\cW)\le k+1$, $L(\cW)\ge L$
and
$\mesh(\cW)\le 6CL$.
This shows that
$\asdim\cone(Z)\le\ldim(Z\times[0,1])$.
\end{proof}

\begin{pro}\label{prop:asprod>} Let
$Z_1,\dots,Z_n$
be bounded metric spaces. Then
$$\ldim(Z_1\times\cdots\times Z_n\times[0,1]^n)
 \le\asdim(\cone Z_1\times\cdots\times\cone Z_n).$$
\end{pro}

We reduce the proof to the following two lemmas.

\begin{lem}\label{lem:asprod>1} Assume that
$\asdim(\cone Z_1\times\cdots\times\cone Z_n)\le k$.
Then for some
$\ov\ep>0$
there exists a function
$\de:[0,\ov\ep]\to [0,1]$
such that for every
$\ep\in(0,\ov\ep)$, $\tau\in(0,1)$
there exists a
$(\de(\ep)\tau,\ep\tau,k+1)$-covering
of
$Z_1\times\cdots\times Z_n\times[0,\tau]^n$
by open subsets of
$Z_1\times\cdots\times Z_n\times\R^n$.
\end{lem}

\begin{proof} We let
$Z=Z_1\times\cdots Z_n$, $X=\cone Z_1\times\cdots\times\cone Z_n$,
and
$h:Z\times\R^n\to X$
the product of the canonical projections. By the assumption,
for every
$L>0$
there exists a
$(L,M,k+1)$-covering
$\cU$
of
$X$
with some
$M<\infty$.
We fix such a covering for
$L>8\de_{\hyp^2}$
together with
$\ep\in(0,\ov\ep)$, $\tau\in(0,1)$
and let
$r=-\ln(\ep\tau)+5M/2+c$,
where the constant
$c=c_0+2\de_{\hyp^2}$
is defined before Claim~\ref{cl:h1}.

The annulus
$\an_{[r,r+M/\ep]}\sub X$
which consists of all
$(x_1,\ldots,x_n)\in X$
with
$r\le|x_io|\le r+M/\ep$, $i=1,\ldots,n$,
is covered by
$\cU(\ep,\tau)=\set{U\in\cU}{$U\cap\an_{[r,r+M/\ep]}\neq\es$}$.
We lift this covering to the covering
$\cU'(\ep, \tau)=h^{-1}(\cU(\ep,\tau))$
of
$Z\times[r,r+M/\ep]^n$
(the later we identify with
$Z\times[0,M/\ep]^n$).

We have
$\mesh(\cU(\ep,\tau))\le M$.
Applying Claim~4(2) with
$l=M$, $R_2=r-M$
(note that
$R_2\ge l$), $R\ge R_2$,
we obtain
$\mesh^\times\cU'(\ep,\tau)\le(\ep\tau,M)$
due to our choices. Applying Claim~4(1)  with
$R_1=r+M/\ep$, $l=L$, $R\le R_1$,
we obtain
$$L^\times\cU'(\ep,\tau)
 \ge(e^{-(r+M/\ep)-L/4-c},L/2)
 =(e^{-M/\ep-L/4-5M/2-2c}\ep\tau,L/2),$$
the later equality holds by the choice of
$r$.

Consider the homothety `along the second factor',
$Z\times[0,M/\ep]^n\to Z\times[0,\tau]^n$,
with coefficient
$\ep\tau/M$
and the image
$\cV$
of
$\cU'$
under this homothety. Then
$\mesh^{\times}(\cV)\le (\ep\tau,\ep\tau)$, $L^\times(\cV)
 \ge(\de_0(\ep)\tau,\ep\tau L/(2M))$
for
$\de_0(\ep)=c_1\ep e^{-M/\ep}$, $c_1=e^{-L/4-5M/2-2c}$,
and we obtain the covering with the required properties taking
$\de(\ep)=\min\{\de_0(\ep),L\ep/(2M)\}$.
\end{proof}

\begin{lem}\label{lem:asprod>2} Let
$Z$
be a metric space. Suppose that for some
$\ov\ep>0$
there exists a function
$\de:[0,\ov\ep]\to [0,1]$
such that for every
$\ep\in(0,\ov\ep)$, $\tau\in(0,1)$
there exists a
$(\de(\ep)\tau,\ep\tau,k+1)$-covering
of
$Z\times[0,\tau]^n$
by open subsets of
$Z\times[0,\infty)^n$.
Then
$\ldim (Z\times[0,1]^n)\le k$.
\end{lem}

\begin{proof}
It suffices to find for sufficiently large
$m\in\N$
a
$(c/m,C/m,k+1)$-covering
of
$Z\times[0,1]^n$
with
$c$, $C$
depending only on
$Z$.

We represent
$$Z\times[0,1]^n=\bigcup_{(l_1,\ldots,l_n)\in A}
  Z\times\left[\frac{l_1}{m},\frac{l_1+1}{m}\right]\times\cdots\times
  \left[\frac{l_n}{m},\frac{l_n+1}{m}\right],$$
where
$A=\{0,\ldots,m-1\}^n$.
We define
$i:A\to S=\{0,\ldots,2^{n}-1\}$
by
$$i(l_1,\ldots,l_n)=(l_1\md 2)2^0+\dots+(l_n\md 2)2^{n-1}.$$
Let
$\ep_0=\min\{1/4,\ov\ep\}$, $\ep_{i+1}=\de(\ep_i)/2$,
$i\in S$.
We fix a
$(\de(\ep_i)/m,\ep_i/m,k+1)$-covering
$\cU_{(l_1,\ldots,l_n)}$
of
$Z\times\left[\frac{l_1}{m},\frac{l_1+1}{m}\right]\times\cdots\times
  \left[\frac{l_n}{m},\frac{l_n+1}{m}\right]$,
where
$i=i(l_1,\ldots,l_n)$.
We put for
$s\in S$
$$\cU_s=\bigcup_{(l_1,\ldots,l_n)\in i^{-1}(s)}\cU_{(l_1,\ldots,l_n)}.$$
In view of the condition
$\ep_{i+1}=\de(\ep_i)/2$
we have
$\mesh(\cU_{s+1})\le L(\cU_s)/2.$
Proceeding by induction over
$s\in S$
and applying Lemma~\ref{lem:union},
we obtain a
$(c/m,C/m,k+1)$-covering
of
$Z\times[0,1]^n$
with
$c=\de(\ep_{2^n-1})$, $C=\ep_0$.
\end{proof}

\begin{proof}[Proof of Proposition~\ref{prop:asprod>}]
Proposition~\ref{prop:asprod>}
follows immediately from previous two lemmas.
\end{proof}

\subsection{Proof of main results}

It is known that any visual hyperbolic space
$X$
can be quasi-isometrically embedded in the hyperbolic cone
over its boundary at infinity
$\di X$
and that for a geodesic hyperbolic space
$Y$
the hyperbolic cone over its boundary at infinity can be
quasi-isometrically embedded in the space itself, see
e.g. \cite{BS3}.

It follows that in these cases
$$\asdim X\le\asdim\cone(\di X)$$
$$\asdim Y\ge\asdim\cone(\di Y)$$
and similar estimates hold for products of such spaces.

Then Theorem~\ref{thm:asprod<} follows immediately from
Proposition~\ref{prop:asprod<}, Theorem~\ref{thm:asprod>}
from Proposition~\ref{prop:asprod>},
Theorem~\ref{thm:lasdim} from
Proposition~\ref{prop:lasdim}.

For the proof of Corollary~\ref{cor:G}, we note that any
Gromov hyperbolic group is a visual hyperbolic space with
the compact, locally self-similar boundary at infinity,
see \cite{BL}. Thus by Proposition~\ref{prop:ssim}, the
spaces of the family
$\set{a\di\Ga_1\times b\di\Ga_2}{$a,b\ge 1$}$
have
$\ldim\le n$
uniformly,
$n=\dim(\di\Ga_1\times\di\Ga_2)$,
for any hyperbolic groups
$\Ga_1$, $\Ga_2$.
By Theorem~\ref{thm:asprod<},
\begin{eqnarray*}
 \asdim(\Ga_1\times\Ga_2)
  &\le&\ldim(\di\Ga_1\times\di\Ga_2)+2\\
  &=&\dim(\di\Ga_1\times\di\Ga_2)+2,
\end{eqnarray*}
where the last equality follows from the fact that
$\di\Ga_1\times\di\Ga_2$
is a compact, locally self-similar space, and
\cite{BL}.
The opposite inequality
$$\asdim(\Ga_1\times\Ga_2)\ge\dim(\di\Ga_1\times\di\Ga_2)+2$$
follows from Theorem~\ref{thm:asprod>} and the inequalities
$$\ldim(Z\times[0,1])\ge\dim(Z\times[0,1])=\dim Z+1,$$
which hold for any metric space
$Z$.

Corollary~\ref{cor:lasdim} follows from Theorem~\ref{thm:asprod>} and
Theorem~\ref{thm:lasdim}.

\section{Applications}

As applications of Corollary~\ref{cor:G}, we have

1) examples of the strict inequality in the product theorem
for the asymptotic dimension in the class of hyperbolic groups.
Namely, it is proved in \cite{Dr1, Dr2} that for every prime
$p$,
there is a hyperbolic Coxeter group
$\Ga_p$
with a Pontryagin surface
$\Pi_p$
as the boundary at infinity. Then by Corollary~\ref{cor:G},
we have

\begin{eqnarray*}
 \asdim(\Ga_p\times\Ga_q)&=&\dim(\Pi_p\times\Pi_q)+2\\
  &<&\dim\Pi_p+\dim\Pi_q+2\\
  &=&\asdim\Ga_p+\asdim\Ga_q
\end{eqnarray*}
for prime
$p\neq q$
(the last equality follows from the main result of \cite{BL});

2) examples of strict inequality in the product theorem
for the hyperbolic dimension (the hyperbolic dimension
is quasi-isometry invariant of metric space introduced
in \cite{BS2}). Indeed, let
$\Ga_p$, $\Ga_q$
for
$p\neq q$
be as in 1). Then

\begin{eqnarray*}
 \hypdim(\Ga_p\times\Ga_q)&\le&\asdim(\Ga_p\times\Ga_q)\\
  &<&\dim\Pi_p+\dim\Pi_q+2\\
  &=&\hypdim(\Ga_p)+\hypdim(\Ga_q)
\end{eqnarray*}
(the last equality also follows from the main result of \cite{BL}
and \cite{BS2});

3) the equality
$\asdim(\Ga\times\R)=\asdim(\Ga)+1$
for any hyperbolic group
$\Ga$.

It is known that for a visual, proper, geodesic hyperbolic  space
$X$
we have
\[\dim(\di X)+1\le\asdim X\le\ldim(\di X)+1 \tag{$\ast$}.\]
The estimate from below is simple and the ideas of the proof are
contained for example in
\cite[1.$\text{E}_1'$]{Gr}.
The estimate from above is proved in \cite{Bu1, Bu2}.
Corollary~\ref{cor:lasdim} gives a better estimate from below and it
allows to show that the first inequality in
($\ast$)
might be strict. More exactly, we give an example
of a hyperbolic space with
the asymptotic dimension arbitrarily larger
than the topological dimension of its boundary at infinity.
Let
$Z=\{0\}\cup\set{1/m}{$m\in\N$}$
and
$X=\cone Z^k$.
It is known that
$\ldim(Z^k\times[0,1])=k+1$,
see \cite{BL}. Then
$\asdim X=k+1$
while
$\dim(\di X)=0$.
The natural conjecture is that for every compact metric space
$Y$
we have
$\ldim(Y\times [0,1])=\ldim(Y)+1$.
If this is true then both inequalities in ($\ast$) become the equalities.

\bigskip
\begin{tabbing}

St. Petersburg Dept. of Steklov\hskip11em\relax \= \\

Math. Institute RAS, Fontanka 27, \> \\

191023 St. Petersburg, Russia\> \\

{\tt lebed@pdmi.ras.ru}\> \\

\end{tabbing}

\end{document}